\documentclass[12pt]{amsart}
\usepackage[margin=0.96in]{geometry}
\usepackage{amsmath}
\usepackage{amssymb}
\usepackage{amsfonts}
\usepackage{xcolor}

\usepackage{mathrsfs}

\newtheorem{theorem}{Theorem}[section]

\newtheorem{lemma}[theorem]{Lemma}

\newtheorem*{conjecture*}{Conjecture}
\newtheorem*{problem*}{Problem}
\newtheorem*{thm*}{Theorem}

\renewcommand{\epsilon}{\varepsilon}

\numberwithin{equation}{section}

\usepackage{mathtools}

\allowdisplaybreaks

\begin{document}

\title{Almost all primes satisfy the Atkin--Serre conjecture and are not extremal}
\author{Ayla Gafni}
\address{Department of Mathematics, The University of Mississippi, Hume Hall 305, University, MS 38677, USA}
\email{ayla.gafni@gmail.com}
\author{Jesse Thorner}
\address{Department of Mathematics, University of Illinois, Urbana, IL 61801, USA}
\email{jesse.thorner@gmail.com}
\author{Peng-Jie Wong}
\address{National Center for Theoretical Sciences, No. 1, Sec. 4, Roosevelt Rd., Taipei City, Taiwan}
\email{pengjie.wong@ncts.tw}


\begin{abstract}
Let $f(z)=\sum_{n=1}^{\infty} a_f(n)e^{2\pi i n z}$ be a non-CM holomorphic cupsidal newform of trivial nebentypus and even integral weight $k\geq 2$.  Deligne's proof of the Weil conjectures shows that $|a_f(p)|\leq 2p^{\frac{k-1}{2}}$ for all primes $p$.  We prove for 100\% of primes $p$ that
\[
2p^{\frac{k-1}{2}}\frac{\log\log p}{\sqrt{\log p}}<|a_f(p)|<\lfloor 2p^{\frac{k-1}{2}}\rfloor.
\]
Our proof gives an effective upper bound for the size of the exceptional set.  The lower bound shows that the Atkin--Serre conjecture is satisfied for 100\% of primes, and the upper bound shows that $|a_f(p)|$ is as large as possible (i.e., $p$ is extremal for $f$) for 0\% of primes.  Our proofs use the effective form of the Sato--Tate conjecture proved by the second author, which relies on the recent proof of the automorphy of the symmetric powers of $f$ due to Newton and Thorne.
\end{abstract}
\maketitle

\section{Introduction and statement of main results}

In this note, we study properties of Fourier coefficients of newforms.  A \emph{newform} $f$ of weight $k$, level $q$, and trivial nebentypus, given by
\begin{equation}
\label{eqn:Fourier_exp}
f(z)=\sum_{n=1}^\infty a_f(n) e^{2\pi i nz}\in S_k^{\mathrm{new}}(\Gamma_0(q)),
\end{equation}
is a normalized cusp form (so that $a_f(1)=1$) that is an eigenform of all the Hecke operators and all of the Atkin--Lehner involutions $|_k W(q)$ and $|_k W(Q_p)$ for each prime $p|q$ (see \cite[Section 2.5]{Ono}).  Throughout, we assume that $f$ is non-CM, i.e., there is no imaginary quadratic field $K$ with the property that $p$ is inert in $K$ if and only if $a_f(p)=0$ (for $p\nmid q$).

The study of the Fourier coefficients $a_f(n)$ of newforms is a central topic in the theory of modular forms.  The motivating case arises from newforms associated to non-CM elliptic curves via modularity.  Indeed, any elliptic curve $E/\mathbb{Q}$ (CM or non-CM) of conductor $q$ has an associated newform
\[
f_E(z) = \sum_{n=1}^{\infty}a_E(n)e^{2\pi i n z}\in S_2^{\mathrm{new}}(\Gamma_0(q)).
\]
This newform encodes information about points on the elliptic curve; for each prime $p\nmid q$, we have the identity $a_E(p) = p+1-\#E(\mathbb{F}_p)$. 
There are also higher-weight newforms which arise naturally, perhaps most notably
\[
\Delta(z) = \sum_{n=1}^{\infty}\tau(n)e^{2\pi i n z}\in S_{12}^{\mathrm{new}}(\Gamma_0(1)),
\]
whose coefficients are given by the Ramanujan tau function.

As a consequence of Deligne's proof of the Weil conjectures, for each prime $p$, there exists $\theta_p\in[0,\pi]$ such that
\[
a_f(p)= 2p^{\frac{k-1}{2}}\cos\theta_p.
\]
It is fruitful to study the distribution of $\cos\theta_p$ as $p$ varies.  The definitive conjecture in this direction is the Sato--Tate conjecture (as extended by Serre \cite{MR1484415}), which was proved by Barnet-Lamb, Geraghty, Harris, and Taylor:  
\begin{thm*}[Sato-Tate Conjecture \cite{BGHT}]  If $H:[-1,1]\to\mathbb{C}$ is Riemann integrable, then
\begin{equation}
\label{eqn:ST_1}
\lim_{x\rightarrow\infty} \frac{1}{\pi(x)}\sum_{p\le x} H(\cos\theta_p)= \int_{-1}^{1}H(t)d\mu_{\mathrm{ST}},\qquad d\mu_{\mathrm{ST}}:=\frac{2}{\pi}\sqrt{1-t^2}~dt,
\end{equation}
where $\pi(x):=\#\{p\leq x\}$ is the usual prime counting function.\footnote{Recall that the prime number theorem states that $\pi(x)\sim x/\log x$.}
\end{thm*}

In a recent breakthrough, Newton and Thorne \cite{NT,NT2} proved that for all integers $n\geq 1$, the $n$-th symmetric power $L$-function $L(s,\mathrm{Sym}^n f)$ associated to $f$ is the $L$-function of an automorphic representation of $\mathrm{GL}_{n+1}(\mathbb{A}_{\mathbb{Q}})$, where $\mathbb{A}_{\mathbb{Q}}$ denotes the ring of adeles over $\mathbb{Q}$.  The zeros of these $L$-functions dictate the distribution of the primes $p$ such that $\theta_p$ lie in a given interval, much like how the zeros of Dirichlet $L$-functions dictate the distribution of primes in arithmetic progressions.  Using the results of Newton and Thorne, the second author in \cite{Th20} proved a strong version of \eqref{eqn:ST_1} with an effective error term.  The purpose of this note is to apply these results to make improvements toward important conjectures concerning the Fourier coefficients $a_f(p)$.  

\subsection{Improvements toward the Atkin-Serre conjecture}
When $k=2$ and $f$ corresponds with a non-CM elliptic curve $E/\mathbb{Q}$, Elkies \cite{Elkies1} proved that $a_f(p)=0$ for infinitely many $p$ (thus $E$ has infinitely many supersingular primes).  In contrast, when $k\geq 4$ we expect that for any fixed $t\in\mathbb{R}$, $a_f(p)=t$ holds for only finitely many primes $p$.  This statement is quantified by a deep conjecture of Atkin and Serre:
\begin{conjecture*}[Atkin-Serre  \cite{Se76}]  Let $f\in S_k^{new} (\Gamma_0(q))$ be a non-CM newform of weight $k\ge 4$.  For each $\epsilon>0$, there exist constants $c_{\epsilon,f}>0$ and $c_{\epsilon,f}'>0$ such that if $p>c_{\epsilon,f}'$, then
\begin{equation}
\label{eqn:AT}
|a_f(p)|\geq c_{\epsilon,f}p^{\frac{k-3}{2}-\epsilon}.
\end{equation}
\end{conjecture*}

Rouse \cite{Rouse} proved for each non-CM newform $f\in S_{k}^{\mathrm{new}}(\Gamma_0(q))$ of weight $k\geq 4$, that \eqref{eqn:AT} holds for almost all primes $p$, conditional on the generalized Riemann hypothesis for the symmetric power $L$-functions associated to $f$.  Unconditionally, M. Ram Murty, V. Kumar Murty, and Saradha \cite{MMS} and V. Kumar Murty \cite{VKM} proved that there exists a density one subset of the primes for which $|a_f(p)|\geq (\log p)^{1-\epsilon}$ for any fixed $\epsilon>0$.  Thus $|a_f(p)|$ cannot be ``too small'' often, but in a much weaker sense than Atkin and Serre predicted.  Our first result establishes \eqref{eqn:AT} unconditionally for almost all primes $p$, a noticeable improvement over the result of Murty, Murty, and Saradha.
\begin{theorem}
	\label{thm:AS_1}
	Let 
	\[
	f(z)=\sum_{n=1}^\infty a_f(n) e^{2\pi i nz}\in S_k^{\mathrm{new}}(\Gamma_0(q))
	\] 
	be a holomorphic normalized cuspidal non-CM newform of even integral weight $k\geq 2$ on $\Gamma_0(q)$ with trivial nebentypus.  
	  There exists an absolute and effectively computable constant $c_1>0$ such that for $x\geq 3$,
	\[
	\#\Big\{x<p\leq 2x\colon |a_f(p)|\leq  2p^{\frac{k-1}{2}}\frac{\log \log p}{\sqrt{\log p}}\Big\}\leq c_1 \frac{x\log(kq \log x)}{(\log x)^{3/2}}.
	\]
	Thus the set of primes $p$ such that \eqref{eqn:AT} holds at $p$ forms a density one subset of the primes.
\end{theorem}

\subsection{Unconditional bounds on extremal primes}
Theorem \ref{thm:AS_1} shows that it is rare for the Fourier coefficients of $f$ to be ``small'' on the primes.  In the opposite direction, one can ask how often the Fourier coefficients of $f$ are ``large'' on the primes.  To make this precise, consider first the holomorphic newforms $f\in S_k^{\mathrm{new}}(\Gamma_0(q))$ such that $a_f(n)\in\mathbb{Z}$ for all $n\geq 1$.  In this case, since $|a_f(p)|\leq 2p^{\frac{k-1}{2}}$, the largest value that $|a_f(p)|$ could assume is $\lfloor 2p^{\frac{k-1}{2}} \rfloor$, where $\lfloor r \rfloor=\max\{t\in\mathbb{Z}\colon t\leq r\}$ is the usual floor function.  If $|a_f(p)|$ assumes this value, we call $p$ {\it extremal} for $f$.  James, Tran, Trinh, Wertheimer, and Zantout \cite{JTTWZ} conjectured that if $f$ is a non-CM newform of weight 2 corresponding with an elliptic curve over $\mathbb{Q}$, then
\begin{equation}
\label{eqn:extremal}
\#\{p\leq x\colon a_f(p) =\lfloor 2\sqrt{p} \rfloor \}\sim \frac{8}{3\pi}\frac{x^{\frac{1}{4}}}{\log x}.
\end{equation}
David, Gafni, Malik, Prabhu, and Turnage-Butterbaugh \cite{Chantal} proved that \eqref{eqn:extremal} is $\ll_f \sqrt{x}$ under the generalized Riemann hypothesis for the symmetric power $L$-functions of $f$.  We prove a nontrivial unconditional upper bound in a more general context than elliptic curves.

\begin{theorem}
	\label{thm:extremal_1}
	Let $f\in S_k^{\mathrm{new}}(\Gamma_0(q))$ be a newform as in Theorem \ref{thm:AS_1}.  There exists an absolute and effectively computable constant $c_3>0$ such that if $x\geq 16$, then
	\[
\#\{x<p\leq 2x\colon |a_f(p)| \geq \lfloor 2p^{\frac{k-1}{2}} \rfloor \}\leq c_3 \frac{x(\log(kq\log x))^2}{(\log x)^2}.
\]
\end{theorem}

\subsection*{Notation}

The Vinogradov notation $F\ll G$ will be used to denote the existence of an effectively computable positive constant $c$ (not necessarily the same in each occurrence) such that $|F|\leq c|G|$ in the range indicated.  We write $F=G+O(H)$ to denote that $|F-G|\ll H$.

\section{Preliminary lemmas}

The main results of this paper follow from the work of Newton and Thorne \cite{NT, NT2} and of the second author \cite{Th20}.  For simplicity, it is assumed that the nebentypus is trivial so that $a_f(n)\in\mathbb{R}$ for all $n\geq 1$; with additional effort, this assumption can be removed.

\begin{lemma}
\label{lem:T1}
	Let $f\in S_k^{\mathrm{new}}(\Gamma_0(q))$ be as in the statement of Theorem \ref{thm:AS_1}, and let $I=[\alpha,\beta]\subseteq[-1,1]$ be a subinterval.  If $x\geq 3$, then
	\[
	\pi_{f,I}(x):=\#\{x<p\leq 2x\colon p\nmid q,~\cos\theta_p\in I\}=(\pi(2x)-\pi(x))\Big(\mu_{\mathrm{ST}}(I)+O\Big(\frac{\log(kq\log x)}{\sqrt{\log x}}\Big)\Big).
	\]
\end{lemma}
	Lemma \ref{lem:T1} is immediate from the work of the second author \cite[Theorem 1.1]{Th20}.  (Note that $\pi(2x)-\pi(x)\sim \pi(x)\sim x/\log x$ by the prime number theorem.)  We will also use the key technical proposition proved in \cite{Th20} which leads to Lemma \ref{lem:T1}.    In what follows, let $U_n$ be the $n$-th Chebyshev polynomial defined by
\[
U_n(\cos\theta) = \frac{\sin((n+1)\theta)}{\sin\theta},\qquad n\geq 0.
\]
The sequence $\{U_n(t)\}_{n=0}^{\infty}$ forms an orthonormal basis of $L^2([-1,1],\mu_{\mathrm{ST}})$.
\begin{lemma}
\label{lem:T2}
	Let $f\in S_k^{\mathrm{new}}(\Gamma_0(q))$ be a newform Lemma \ref{lem:T1}.  There exist absolute and effectively computable constants $c_5>0$ and $c_6>0$ such that if $1\leq n\ll \sqrt{\log x}/\sqrt{\log(kq\log x)}$, then
	\[
	\Big|\sum_{x<p\leq 2x}U_n(\cos\theta_p)\Big|\ll \frac{x}{\log x}\Big(n^2 x^{-\frac{1}{c_5 n}}+n^2\Big(\exp\Big[-c_6\frac{\log x}{n^2\log(kqn)}\Big]+\exp\Big[-c_6\frac{\sqrt{\log x}}{\sqrt{n}}\Big]\Big)\Big).
	\]
\end{lemma}
	Lemma \ref{lem:T2} follows immediately from \cite[Proposition 2.1]{Th20} and partial summation.
	
\section{Proof of Theorem \ref{thm:AS_1}}

Since $a_f(p) =   2p^{\frac{k-1}{2}}\cos\theta_p$, we see that
\begin{align*}
	\#\Big\{x<p\leq 2x\colon |a_f(p)|\leq 2p^{\frac{k-1}{2}}\frac{\log(kq\log p)}{\sqrt{\log p}}\Big\}\leq \#\Big\{x<p\leq 2x\colon |\cos\theta_p|\leq\frac{\log(kq\log x)}{\sqrt{\log x}}\Big\}.
\end{align*}
Writing $I=[-\frac{\log(kq\log x)}{\sqrt{\log x}},\frac{\log(kq\log x)}{\sqrt{\log x}}]$, we find via a Taylor expansion that
\[
\mu_{\mathrm{ST}}(I)=\frac{4}{\pi}\frac{\log(kq\log x)}{\sqrt{\log x}}+O\Big(\frac{(\log(kq\log x))^3}{(\log x)^{3/2}}\Big).
\]
Thus the desired result follows from Lemma \ref{lem:T1}.

\section{Proof of Theorem \ref{thm:extremal_1}}\label{proof-main-3}

For an integer $N\geq 3$ (which we will later take to depend on $x$), let $I_{N}=[\cos N^{-1},1]$ and $I_{N}'=[-1,-\cos N^{-1}]$.  We take $N$ large enough so that $\cos N^{-1}\leq 1-x^{-\frac{1}{2}}$.  If $x<p\leq 2x$, then we have $x< p<4p^{k-1}$.  With $\{t\}=t-\lfloor t\rfloor$, we deduce the chain of inequalities
\[
\cos N^{-1}\leq 1-x^{-\frac{1}{2}}<1-\frac{1}{2p^{\frac{k-1}{2}}}<1-\frac{\{2p^{\frac{k-1}{2}}\}}{2p^{\frac{k-1}{2}}}<1.
\]
Therefore, $\#\{x< p \le 2x\colon |a_f(p)|= \lfloor 2p^{\frac{k-1}{2}}\rfloor \}$ equals
\begin{multline*}
\#\Big\{ x< p \le 2x\colon  \frac{|a_f(p)|}{2p^{\frac{k-1}{2}}}
=  1-\frac{\{2p^{\frac{k-1}{2}}\}}{2p^{\frac{k-1}{2}}} \Big\}\\
\le \#\{x< p \le 2x\colon \cos\theta_p \in I_{N}\}+\#\{x< p \le 2x\colon \cos\theta_p \in I_{N}'\}.
\end{multline*}
We will estimate $\#\{x< p \le 2x\colon \cos\theta_p \in I_{N}\}$; the case for $I_{N}'$ is the same.  By the work in \cite{Chantal} (Proposition 2.2, Equation 3.2, and the displayed equation preceding it), we find that if $N\geq 3$, then (note that $U_0(\cos\theta_p)$ is identically $1$)
\begin{equation}
\label{eqn:last_step}
\#\{x< p \le 2x\colon  \cos\theta_p \in I_{N}\}\ll \frac{\pi(2x)-\pi(x)}{N^2}+\frac{1}{N^2}\sum_{n=1}^{N}\Big|\sum_{x<p\leq 2x}U_n(\cos\theta_p)\Big|.
\end{equation}
Theorem \ref{thm:extremal_1} follows once we apply Lemma \ref{lem:T2} to \eqref{eqn:last_step} with $N=\lfloor\sqrt{\log x}/\log(kq\log x)\rfloor$.

\section*{Acknowledgments}
The authors would like to thank Amir Akbary, Po-Han Hsu, and Wen-Ching Winnie Li for helpful comments.  The third author is currently an NCTS postdoctoral fellow; he was supported by a PIMS postdoctoral fellowship and the University of Lethbridge during part of this research.

\bibliographystyle{abbrv}
\bibliography{Bibliography}

\end{document}